\documentclass[12pt,reqno,final,a4paper]{amsart}

\usepackage[T1]{fontenc}
\usepackage{lmodern}
\usepackage{microtype}
\usepackage[margin=1in,footskip=24pt]{geometry}

\usepackage{xcolor}
\definecolor{darkgray}{gray}{0.9}

\usepackage{amsmath,amssymb,amsthm}
\usepackage{mathtools,mathrsfs,bbm}
\usepackage{thmtools,thm-restate}

\allowdisplaybreaks

\usepackage{graphicx}
\usepackage{tikz}
\usetikzlibrary{fit,backgrounds,shapes,calc,arrows.meta}
\usepackage{tikzscale}
\usepackage{pgfplots}
\pgfplotsset{compat=1.12}

\usepackage{pifont}

\usepackage[linesnumbered,ruled,vlined]{algorithm2e}
\usepackage[noend]{algpseudocode}

\SetCommentSty{mycommfont}
\SetKwInOut{Input}{input}
\SetKwInOut{Output}{output}

\usepackage[shortlabels]{enumitem}
\usepackage{environ}
\usepackage{framed}
\usepackage[framemethod=tikz]{mdframed}
\usepackage{tcolorbox}
\usepackage[labelfont=bf]{caption}
\usepackage{appendix}
\usepackage{textcomp}
\usepackage{ifthen}
\usepackage[textsize=tiny]{todonotes}

\usepackage[noadjust]{cite}
\usepackage[hidelinks,hypertexnames=false]{hyperref}
\usepackage[noabbrev,capitalize,nameinlink]{cleveref}

\crefname{equation}{}{}

\crefformat{enumi}{#2#1#3}
\crefrangeformat{enumi}{#3#1#4 to~#5#2#6}
\crefmultiformat{enumi}{#2#1#3}
    { and~#2#1#3}
    {, #2#1#3}
    { and~#2#1#3}

\makeatletter
\newcommand{\namedlabel}[2]{%
    \phantomsection
    \def\@currentlabelname{#2}%
    \label{#1}%
}
\makeatother

\numberwithin{equation}{section}

\theoremstyle{plain}
\newtheorem{theorem}{Theorem}[section]

\newtheorem{lemma}[theorem]{Lemma}
\newtheorem{claim}[theorem]{Claim}

\newtheorem{conjecture}[theorem]{Conjecture}
\newtheorem*{question*}{Question}

\theoremstyle{definition}

\newtheorem*{definition*}{Definition}

\newtheorem*{fact*}{Fact}

\theoremstyle{remark}

\crefname{claim}{Claim}{Claims}
\crefname{fact}{Fact}{Facts}

\DeclareMathOperator{\Aut}{Aut}

\DeclareMathOperator{\Bin}{Bin}
\DeclareMathOperator{\Var}{Var}

\DeclareMathOperator{\deff}{def}

\newcommand{\tiso}{\tau_{\mathrm{iso}}}

\newlength{\entropynamewidth}

\title{The number of non-isomorphic spanning trees in regular graphs}
\author{Csongor Beke}
\date{}
\thanks{University of Cambridge, \texttt{cb2138@cam.ac.uk}.\quad \today.}

\begin{document}

\begin{abstract}
By a classical theorem of Otter, the number of non-isomorphic trees on $n$ vertices is $\alpha^n\exp(-o(n))$, where $\alpha\approx 2.956$ is Otter's constant. We prove that every connected $d$-regular graph $G$ on $n$ vertices contains $\alpha^n\exp(-o_d(n))$ non-isomorphic spanning trees, where $o_d(n)/n\to0$ as $d\to\infty$.
\end{abstract}

\maketitle

\section{Introduction}

How many spanning trees does a connected graph $G$ contain? Write $\tau(G)$ for this quantity, counting trees with different edge sets separately. Cayley's formula~\cite{Cayley} gives $\tau(K_n)=n^{n-2}$, while Kirchhoff's matrix-tree theorem~\cite{kirchhoff1847ueber} expresses $\tau(G)$ in terms of the Laplacian of $G$. For connected $d$-regular graphs on $n$ vertices, the work of McKay~\cite{McKay} and Alon~\cite{Alon} gives $\tau(G)=d^n\exp(-o_d(n))$. Kostochka~\cite{kostochka1995number} extended these estimates to graphs with a prescribed degree sequence. Throughout the paper, all graphs are simple, and we take $d\geq1$.

We study the corresponding question when spanning trees are counted up to isomorphism. Write $\tiso(G)$ for the number of isomorphism classes of spanning trees of $G$. For the complete graph, this is precisely the number of unlabelled trees on $n$ vertices. Otter's theorem gives asymptotic formulae for this number and for the number $r_n$ of non-isomorphic rooted trees on $n$ vertices.
\begin{theorem}[Otter~\cite{Otter}]\label{thm:otter}
There exist absolute constants $\alpha,C,C_{\mathrm r}>0$ so that, as $n\to\infty$,
\[
    \tiso(K_n)=(1+o(1))Cn^{-5/2}\alpha^n,
    \qquad
    r_n=(1+o(1))C_{\mathrm r}n^{-3/2}\alpha^n.
\]
\end{theorem}

The study of $\tiso(G)$ beyond complete graphs goes back to the 1970s, when Zelinka~\cite{Zelinka1971,Zelinka1978} and  Hartnell~\cite{Hartnell1974} investigated graphs with few isomorphism classes of spanning trees. Recently, Lee~\cite{Lee} proved that every connected almost regular graph of sufficiently large degree contains exponentially many non-isomorphic spanning trees. He asked whether, for connected $d$-regular graphs, the exponential base tends to $\alpha$ as $d\to\infty$. Our main result answers this question affirmatively.
\begin{theorem}\label{thm:main}
Every connected $d$-regular graph $G$ on $n$ vertices satisfies
\[
    \tiso(G)= \alpha^n\exp\bigl(-o_d(n)\bigr)
    =\bigl(\alpha-o_d(1)\bigr)^n,
\]
where the error tends to zero as $d\to\infty$, uniformly in $n$ and $G$.
\end{theorem}
For fixed $d$, an exponential loss in $n$ is unavoidable; see \cref{sec:remarks} for more details.

One approach to counting non-isomorphic spanning trees in a $d$-regular graph $G$ is to bound, for each fixed $n$-vertex tree $T$, the number of spanning trees of $G$ that are isomorphic to $T$. Lee~\cite{Lee} obtained the bound $d^n\exp(-\Omega(n))$, while Beke, Bo\v{s}kovi\'c, Kam\v{c}ev and Wang~\cite{BBKW} subsequently improved this to $d^n\exp(-n+o_d(n))$. Combined with $\tau(G)=d^n\exp(-o_d(n))$, these bounds give $\tiso(G)\geq e^{\Omega(n)}$ and $\tiso(G)\geq e^{n-o_d(n)}$, respectively. The bound on the size of a single isomorphism class is sharp up to the $o_d(n)$ term, even for paths, so this argument cannot yield an exponential base larger than $e<\alpha$.

The conclusion of \cref{thm:main} fails if regularity is replaced by a minimum-degree assumption: for fixed $d$, the complete bipartite graph $K_{d,n-d}$ has only polynomially many spanning-tree classes as $n\to\infty$. In this setting, Bitonti, Michel and Scott~\cite{BitontiMichelScott} proved a lower bound of $n^{\Omega(d)}$ for connected $n$-vertex graphs of minimum degree at least $d$, when $d$ is sufficiently large and $n$ is sufficiently large relative to $d$. Yan and Zhang~\cite{Yan-Zhang} subsequently showed that, for every fixed $d\geq3$ and sufficiently large $n$, the unique minimiser is $K_{d,n-d}$. The problem of determining $\tiso(K_{a,b})$ has also been studied, with exact formulae in special cases and bounds for general $a$ and $b$;
see~\cite{ClarkMohrPorter2007,vanDenBoomen2009,JohnsonNochumson2026}.

A different approach for regular graphs was taken by Wang~\cite{Wang}, who proved the conclusion of \cref{thm:main} under an additional spectral assumption. His proof constructs many pairwise non-isomorphic trees on $n$ vertices by joining bounded-degree rooted trees along a path, and then uses a spanning-tree universality theorem to embed them into pseudorandom graphs. Our proof similarly joins rooted trees along a path-like structure, but constructs them directly inside the given $d$-regular graph, avoiding the need for a universality theorem. Different assignments of rooted trees may produce isomorphic spanning trees, and the main technical challenge is to bound the number of assignments yielding any one isomorphism class.
\section{Proof sketch}
Our proof constructs spanning trees by joining many rooted trees of size $k=\lfloor\sqrt{\log d}\rfloor$. We use the rooted part of  \cref{thm:otter}, which gives $r_k=\alpha^k\exp(-o_k(k))$. Let $Q_1,\ldots,Q_{r_k}$ be representatives of the rooted tree classes on $k$ vertices. Label the vertices of each tree by $0,\ldots,k-1$, with the root labelled $0$.

\subsection{The complete graph case} To demonstrate our method, we first show how to obtain $\tiso(K_n)\geq \alpha^n\exp(-o_k(n))$ from  $r_k= \alpha^k\exp(-o_k(k))$. Take a path $S$ on $n/k$ vertices, and for each vertex of $S$, attach a copy of one of the rooted trees $Q_s$, identifying its root with the vertex of $S$. These assignments give us $r_k^{n/k}$ spanning trees of $K_n$, but some of them could be isomorphic. To bound how many assignments can produce some given tree $T^{\ast}$, we first specify the endpoints of $S$ and their order inside $T^{\ast}$. This can be done in at most $2\binom{n}{2}$ ways, and after we have identified $S$ in $T^{\ast}$, the assignments of the rooted trees can be read off: delete the edges of $S$; each resulting component, rooted at its unique vertex of $S$, reveals the assigned rooted tree class. Therefore, we obtain
\[
    \tiso(K_n)\geq \frac{r_k^{n/k}}{2\binom{n}{2}}=\alpha^n\exp\bigl(-o_k(n)\bigr).
\]

\subsection{Packing vertex-disjoint rooted trees} For a general connected $d$-regular graph, we would like to follow a similar strategy; the difficulty is embedding the trees. We start by finding these disjoint copies of rooted trees in $G$, and then complete them to spanning trees.

Our approach of finding the copies of $Q_s$ is similar to the method of Kelmans, Mubayi and Sudakov~\cite{KMS} for packing regular graphs with copies of a fixed tree. We consider an ordered partition of $V(G)$ into $kr_k$ pairwise disjoint parts,
\[
    \Pi=(V_{s,j}:s\in[r_k],\ 0\leq j\leq k-1),
\]
where we think of $V_{s,j}$ as the set of possible images of vertex $j$ in copies of $Q_s$. Then, for each $s$ and each edge $ij\in E(Q_s)$ we take a maximum matching between $V_{s,i}$ and $V_{s,j}$. Their union $H_{\Pi}$ is a spanning forest of $G$, whose components are isomorphic to subtrees of the corresponding rooted trees. We call a component \emph{complete} if it has $k$ vertices. Let $R=\bigcup_s V_{s,0}$ be the root set, and assign $Q_s$ to each root in $V_{s,0}$. Every complete component contains exactly one root and is a copy of the tree assigned to it, while each incomplete component contains at most one root.

Crucially, writing $D_{\Pi}$ for the total \emph{deficiency}---the sum of the numbers of vertices left unmatched by each matching---at most $D_\Pi$ components are incomplete. For a uniformly random partition $\Pi$, the expected deficiency is small; see \cref{lem:matching}.

\subsection{Constructing many spanning trees}
Next, we need to extend $H_\Pi$ to a spanning tree of $G$. Since $G$ is connected, this is always possible, but we also want to recover the assigned rooted trees from the resulting tree. In the complete graph case, a path joins the roots. Here we use a path forest.

For now, fix the root set $R$. Choose a path forest $F_R$ in $G[R]$ covering $R$ with the minimum possible number $c(R)$ of components. Since each component of $H_\Pi$ contains at most one root, $F_R\cup H_\Pi$ is a forest. We extend it to a spanning tree $T_\Pi$ by adding at most $D_\Pi+c(R)$ edges.

Call the minimal subtree of $T_\Pi$ containing $R$ the \emph{backbone}, and denote it by $S_\Pi$. Since $S_\Pi$ contains the fixed forest $F_R$, we can specify it inside $G$ by recording its edges outside $F_R$. As in the complete graph case, deleting the backbone edges recovers the assigned rooted tree at every root whose component of $H_\Pi$ is complete and meets no added edge.

\subsection{Good partitions}Set $h=n/k^4$, and call $\Pi$ \emph{good} if $D_\Pi+c(R)\leq h$. For a good partition, at most $h$ components of $H_\Pi$ are incomplete, and at most $h$ edges are added when constructing $T_\Pi$. These edges meet at most $2h$ further components, so all but at most $3h$ components remain complete copies attached to the backbone only at their roots.

By considering uniformly random partitions, we show that at least half of all partitions are good; see \cref{claim:good}. An averaging argument then gives a root set $R$ for which at least half of all assignments of rooted tree classes to $R$ extend to a good partition of $V(G)$. Fix one such extension for each of these assignments. For this choice of $R$, we have $|R|=n/k+O(h)$, so we retain at least $r_k^{n/k}\exp(-o_d(n))$ assignments.

\subsection{Recovering the assignments}We record four pieces of information that allow us to recover the rooted class assignment from a fixed tree $T^\ast$. First, we identify $S_\Pi$ inside $G$. Since it contains $F_R$ and at most $3kh$ further edges, there are only $\exp(o_d(n))$ possibilities. Second, we specify the image of $S_\Pi$ inside $T^\ast$. The backbone has at most $2h$ leaves and is determined by their images, giving another $\exp(o_d(n))$ possibilities. Third, we specify an isomorphism between these two copies of the backbone. There are at most $|\Aut(S_\Pi)|$ choices, which we bound using \cref{lem:aut}. Finally, deleting the edges of $S_\Pi$ reveals the assignment at each root whose component of $H_\Pi$ is complete and meets no added edge. There are at most $3h$ other roots, so we can specify them and their assignments in at most $\sum_{i\le 3h}\binom{|R|}{i}r_k^i=\exp(o_d(n))$ ways.

Thus the multiplicity is at most $\exp(o_d(n))$, and dividing by this bound gives 
\[
    \tiso(G)\geq r_k^{n/k}\exp(-o_d(n))=\alpha^n\exp(-o_d(n)).
\]

\section{Auxiliary lemmas}
    \label{sec:auxiliary}
We first bound the expected number of vertices left unmatched by a maximum matching between two random vertex sets in a regular graph.

For disjoint vertex sets $U,V\subset V(G)$, let $\nu(U,V)$ be the size of a maximum matching in the induced bipartite graph between $U$ and $V$. Write
\[
    \deff(U,V)=|U|+|V|-2\nu(U,V)
\]
for its \emph{deficiency}, that is, the number of unmatched vertices in a maximum matching.

\begin{lemma}\label{lem:matching}
Let $G$ be a $d$-regular graph on $n$ vertices, and let $p\in(0,1/2]$. Form disjoint random sets $U,V\subseteq V(G)$ by assigning each vertex independently to $U$, to $V$, or to neither, with probabilities $p$, $p$ and $1-2p$, respectively. Then
\[
    \mathbb E\deff(U,V)\leq 2n\sqrt{p/d}.
\]
\end{lemma}

\begin{proof}
    Let $H$ be the bipartite graph induced by $G$ between $U$ and $V$. By K\"onig's theorem, there is a vertex cover $C\subseteq U\cup V$ of size $\nu(U,V)$. Then $I=(U\cup V)\setminus C$ is independent, so $\sum_{v\in I}d_H(v)\leq\sum_{v\in C}d_H(v)$. Therefore, as $\deff(U,V)=|I|-|C|$, 
\[
    pd\deff(U,V)=pd(|I|-|C|)\leq \sum_{v\in I}(pd-d_H(v))+\sum_{v\in C}(d_H(v)-pd)\leq \sum_{v\in U\cup V}|d_H(v)-pd|.
\]
Conditional on $v\in U\cup V$, the degree $d_H(v)$ is distributed as $\Bin(d,p)$.  By Cauchy--Schwarz,
\[
    \mathbb E|\Bin(d,p)-dp|
    \leq\sqrt{\Var(\Bin(d,p))}\leq\sqrt{dp}.
\]
Taking expectations gives
\[
    \mathbb E\deff(U,V)
    \leq\frac{2np}{dp}\mathbb E|\Bin(d,p)-dp|
    \leq2n\sqrt{p/d}.\qedhere
\]
\end{proof}

We next use \cref{lem:matching} to find a path forest with few components covering a random vertex set. For $R\subseteq V(G)$, write $c(R)$ for the minimum number of vertex-disjoint paths in $G[R]$ covering $R$, allowing single vertices as paths.

\begin{lemma}\label{lem:paths}
Let $G$ be a $d$-regular graph on $n$ vertices, and let $p\in(0,1]$. Form a random set $R\subseteq V(G)$ by including each vertex independently with probability $p$. Then
\[
    \mathbb E c(R)\leq \frac{3n}{d^{1/3}}.
\]
\end{lemma}

\begin{proof}
Put $t=\lfloor d^{1/3}\rfloor$ and partition $R$ into sets $U_1,\ldots,U_t$ by assigning each vertex independently and uniformly to one of these sets. For each $i\in[t-1]$, take a maximum matching $M_i$ between $U_i$ and $U_{i+1}$. Their union is a spanning path forest on $R$, so
\[
    c(R)\leq |R|-\sum_{i=1}^{t-1}|M_i|
    =\frac{|U_1|+|U_t|}{2}
     +\frac12\sum_{i=1}^{t-1}\deff(U_i,U_{i+1}).
\]
Each vertex of $G$ lies in a given $U_i$ with probability $p/t$, independently of the other vertices. Applying \cref{lem:matching} to each consecutive pair and using $d^{1/3}/2\leq t\leq d^{1/3}$, we get
\[
    \mathbb E c(R)
    \leq \frac{np}{t}+n(t-1)\sqrt{\frac{p}{td}}
    \leq \frac nt+n\sqrt{\frac td}
    \leq \frac{3n}{d^{1/3}}.\qedhere
\]
\end{proof}
 The following standard bound will control the number of automorphisms of the backbone.

\begin{lemma}\label{lem:aut}
For $d\geq1$ and $\ell\geq2$, every tree $S$ with maximum degree at most $d$ and at most $\ell$ leaves satisfies
\[
    \lvert\Aut(S)\rvert\leq\ell d^{\ell-2}.
\]
\end{lemma}

\begin{proof}If $S$ has fewer than two vertices the result is clear. Otherwise, fix a leaf $r$ and root $S$ at $r$. Under an automorphism, $r$ can be mapped to at most $\ell$ vertices. After choosing its image, we choose the images of the other vertices one by one, away from the root.

The image of the neighbour of $r$ is determined, and whenever we reach a vertex of degree two, the image of its child is determined as well. At a vertex of degree $q\geq3$, we have at most $d$ choices for the image of each of its first $q-2$ children, after which the image of the last child is determined. Since
\[
    \sum_{v:\deg_S(v)\geq3}\bigl(\deg_S(v)-2\bigr)\leq\ell-2,
\]
we have at most $\ell d^{\ell-2}$ automorphisms.\end{proof}

\section{Proof of Theorem~\ref{thm:main}}
Let $G$ be a connected $d$-regular graph on $n$ vertices. We may assume that $d$ is sufficiently large, and set
\[
    k=\lfloor\sqrt{\log d}\rfloor,     \qquad h=\frac{n}{k^4}.
\]
Let $Q_1,\ldots,Q_{r_k}$ be a list of all non-isomorphic rooted trees on $k$ vertices. Label the vertices of each tree by $0,\ldots,k-1$, with the root labelled $0$. 

Let $\mathcal P$ be the family of all ordered partitions of $V(G)$ of the form
\[
    \Pi=(V_{s,j}:s\in[r_k],\ 0\leq j\leq k-1),
\]
allowing empty parts. We think of $V_{s,j}$ as the possible images of vertex $j$ in copies of $Q_s$. For each $s\in[r_k]$ and $ij\in E(Q_s)$, choose a maximum matching between $V_{s,i}$ and $V_{s,j}$, and let $H_\Pi$ be the spanning graph on $V(G)$ consisting of the union of these matchings. Fix these choices for every partition and set
\[
    D_\Pi=\sum_{s\in[r_k]}\sum_{ij\in E(Q_s)}
           \deff(V_{s,i},V_{s,j}),
    \qquad R(\Pi)=\bigcup_{s\in[r_k]}V_{s,0}.
\]

The graph $H_\Pi$ is a forest whose components are copies of subtrees of the $Q_s$, so they have at most $k$ vertices. Call a component \emph{complete} if it has $k$ vertices; it is then a copy of its assigned rooted tree. Every incomplete component contains a vertex unmatched in at least one prescribed matching. Thus there are at most $D_\Pi$ incomplete components, containing at most $kD_\Pi$ vertices in total.

\subsection{Good partitions} For each $R\subseteq V(G)$, recall that $c(R)$ is the minimum number of vertex-disjoint paths in $G[R]$ covering $R$. Define a partition $\Pi$ to be \emph{good} if
\[
    D_\Pi+c(R(\Pi))\le h.
\]
Our choice of $h$ is small enough to control the resulting trees, but large enough that at least half of all partitions are good.

\begin{claim}\label{claim:good}
For $d$ sufficiently large, at least half of the partitions in $\mathcal P$ are good.
\end{claim}

\begin{proof}
Recall that $\log r_k=O(k)$, so our choice of $k$ gives $kr_k\leq d^{1/3}$ for sufficiently large $d$. Choose $\Pi$ uniformly from $\mathcal P$, or equivalently assign each vertex independently and uniformly to one of the $kr_k$ parts. Applying \cref{lem:matching} with $p=1/(kr_k)$ to each prescribed pair gives
\[
    \mathbb E D_\Pi \leq\frac{2nr_k(k-1)}{\sqrt{kr_kd}}\leq2n\sqrt{\frac{kr_k}{d}}\leq2nd^{-1/3}.
\]
Each vertex belongs to $R(\Pi)$ independently with probability $1/k$. By \cref{lem:paths}, $\mathbb E c(R(\Pi))\leq 3nd^{-1/3}$. Hence
\[
    \mathbb E\bigl[D_\Pi+c(R(\Pi))\bigr]\leq5nd^{-1/3}=o(h),
\]
and Markov's inequality finishes the proof.
\end{proof}
For $R\subseteq V(G)$, let $\Sigma(R)$ be the set of functions $\sigma:R\to[r_k]$ that extend to a good partition $\Pi$, with $V_{s,0}=\sigma^{-1}(s)$ for each $s\in[r_k]$.

By \cref{claim:good} and averaging over root sets, there is some $R$ for which at least half of the partitions with root set $R$ are good. Each function $\sigma:R\to[r_k]$ has the same number of extensions to a partition with root set $R$, so
\begin{equation}\label{eq:sigma}
    |\Sigma(R)|\geq r_k^{|R|}/2.
\end{equation}
\subsection{Constructing spanning trees} Fix such an $R$, and let $F$ be a path forest in $G[R]$ covering $R$ with $c(R)$ components. For each $\sigma\in\Sigma(R)$, choose one good extension $\Pi_\sigma$. Write $H_\sigma=H_{\Pi_\sigma}$ and $D_\sigma=D_{\Pi_\sigma}$, so that
\[
c(R)+D_\sigma\leq h.
\]

Each component of $H_\sigma$ contains at most one root, so $F\cup H_\sigma$ is a forest. Since $G$ is connected, we can extend it to a spanning tree $T_\sigma$ of $G$. Fix such a choice for each $\sigma$, and let $E_\sigma$ be the set of added edges.

We have now constructed a family of spanning trees indexed by $\Sigma(R)$, of size at least $r_k^{|R|}/2$. It remains to bound $|R|$ and the number of $\sigma\in \Sigma(R)$ producing the same isomorphism class. To bound this, we will recover $\sigma$ from the resulting tree and a small amount of additional information. We first establish the required bounds.

For the size of $R$, note that each complete component of $H_\sigma$ has $k$ vertices and exactly one root. There are at most $D_\sigma\le h$ incomplete components, each with at most $k$ vertices and at most one root. Therefore,
\begin{equation}\label{eq:R-size}
    \frac nk-h\leq|R|\leq\frac nk+h\leq\frac{2n}{k}.
\end{equation}

The components of $H_\sigma$ containing roots are joined by $F$ into exactly $c(R)$ components. Every remaining component is rootless and therefore incomplete, so there are at most $D_\sigma$ of them. Thus $|E_\sigma|\leq c(R)+D_\sigma-1\leq h$.

Let $S_\sigma$ be the minimal subtree of $T_\sigma$ containing $R$, which we call the backbone. It contains $F$, and every leaf of $S_\sigma$ lies in $R$ and cannot be an internal vertex of a path of $F$. Thus $S_\sigma$ has at most $2c(R)\leq2h$ leaves.

Call a component of $H_\sigma$ \emph{exceptional} if it is incomplete or contains an endpoint of an edge of $E_\sigma$. There are at most $D_\sigma+2|E_\sigma|\leq3h$ exceptional components. Let $B_\sigma$ be the set of roots lying in these components, so $|B_\sigma|\leq3h$.

Every non-exceptional component is a complete copy of its assigned rooted tree. Since no edge of $E_\sigma$ meets it, it is joined to the rest of $T_\sigma$ only through its root, and hence meets $S_\sigma$ only at that root.

Every edge of $S_\sigma$ outside $F$ lies either in $E_\sigma$ or in an exceptional component of $H_\sigma$. Since there are at most $3h$ exceptional components, each with at most $k-1$ edges, we get
\begin{equation}\label{eq:C-size}
    |E(S_\sigma)\setminus E(F)|
    \leq3h(k-1)+|E_\sigma|\leq3kh.
\end{equation}

\subsection{Recovering the assignments} Fix an $n$-vertex tree $T^\ast$. For every $\sigma\in\Sigma(R)$ with $T_\sigma\cong T^\ast$, choose an isomorphism to $T^\ast$ and make
the following four records.
\begin{itemize}
    \item First, record the edges of $S_\sigma$ outside $F$ as a subset of $E(G)$. By \eqref{eq:C-size}, there are at most $\sum_{i\leq3kh}\binom{nd/2}{i}$ choices. Since $F$ is fixed in $G$, this determines the backbone $S_\sigma$ as a subtree of $G$.
    \item Second, record the images of the leaves of $S_\sigma$ in $T^\ast$. There are at most $\sum_{i\leq2h}\binom ni$ choices. A subtree of a tree is determined by its leaves, so this identifies the image of the backbone in $T^\ast$.
    \item Third, record the isomorphism between these two copies of the backbone. There are $|\Aut(S_\sigma)|$ choices. Since $S_\sigma$ has maximum degree at most $d$ and at most $2h$ leaves, \cref{lem:aut} bounds this by $2h\,d^{2h-2}$.
    \item Finally, record $B_\sigma$ and the values of $\sigma$ on $B_\sigma$. There are at most $\sum_{i\leq3h}\binom{|R|}{i}r_k^i$ choices.
\end{itemize}

For every $r\in R\setminus B_\sigma$, deleting the edges of the backbone from $T^\ast$ leaves a component rooted at $r$ that is isomorphic to $Q_{\sigma(r)}$. Thus these records determine $\sigma$. Consequently, at most $M$ assignments in $\Sigma(R)$ produce any given isomorphism class, where
\[
    M=
    \left(\sum_{i\leq3kh}\binom{nd/2}{i}\right)
    \left(\sum_{i\leq2h}\binom ni\right)
    2h\,d^{2h-2}
    \left(\sum_{i\leq3h}\binom{|R|}{i}r_k^i\right).
\]

To bound $\log M$, we apply the following bound to each sum:
\[
    \sum_{i\le x}\binom Ni\le(eN/x)^x
    \qquad\text{for }0<x\le N.
\]
Together with $h=n/k^4$ and $\log d=\Theta(k^2)$, the first factor gives
\[
    \log\left(\sum_{i\leq3kh}\binom{nd/2}{i}\right)
    \leq3kh\log\left(\frac{end}{6kh}\right)
    =O\left(\frac{n}{k^3}\log(dk^3)\right)
    =O(n/k).
\]
The next two factors give
\[
    \log\left[\left(\sum_{i\leq2h}\binom ni\right)
                   2h\,d^{2h-2}\right]
    =O\bigl(h\log(n/h)+h\log d\bigr)
    =O(n/k^2).
\]
Finally, $|R|\leq2n/k$ and $\log r_k=O(k)$ imply
\[
    \log\left(\sum_{i\leq3h}\binom{|R|}{i}r_k^i\right)
    \leq3h\log\left(\frac{er_k|R|}{3h}\right)
    =O\bigl(h\log(r_kk^3)\bigr)
    =O(n/k^3).
\]
Hence $\log M=O(n/k)$. Combining this with
\eqref{eq:R-size} and~\eqref{eq:sigma}, we obtain
\[
    \log\tiso(G)\geq\log|\Sigma(R)|-\log M
    \geq\frac nk\log r_k-O(n/k)
    =n\log\alpha-o_d(n).
\]
This proves the lower bound, uniformly in $n$ and $G$. The matching upper bound follows from $\tiso(G)\leq\tiso(K_n)$ and \cref{thm:otter}, since $n\geq d+1$.

\section{Concluding remarks}
\label{sec:remarks}
For fixed $d\geq3$, the number of unlabelled trees on $n$ vertices with maximum degree at most $d$ is $\alpha_d^{n}\exp(o(n))$, where $\alpha_d<\alpha$ and $\log(\alpha/\alpha_d)=\Theta(\alpha^{-d})$ as $d\to\infty$~\cite{GohSchmutz1994}. Since every spanning tree of a $d$-regular graph has maximum degree at most $d$, this already shows that an exponential loss in $n$ is unavoidable in \cref{thm:main} for fixed $d$. The following construction shows that this loss can be much larger.

\begin{lemma}\label{prop:necessary-error}
For every sufficiently large $d$ and arbitrarily large $n$, there is a connected $d$-regular graph $G$ on $n$ vertices with $\tiso(G)\leq\alpha^n\exp\bigl(-(\tfrac52+o_d(1))n\log d/d\bigr)$.
\end{lemma}

\begin{proof}
Put $k=d+3$ and form $B$ from $K_k$ by deleting the edges of a triangle and a $(k-2)$-cycle meeting in a single vertex $x$. Then $x$ has degree $d-2$ and all other vertices have degree $d$. Take a cycle $u_1v_1\cdots u_mv_mu_1$ and, for each $i$, join both $u_i$ and $v_i$ to $x$ in each of $d-2$ fresh copies of $B$. The resulting graph $G$ is connected and $d$-regular, with $n=m(k(d-2)+2)$ vertices.

Since only $x$ has neighbours outside each copy of $B$, every spanning tree of $G$ restricts to a spanning tree of each copy. For each copy, record the rooted class of the restriction at $x$ and which of $xu_i,xv_i$ belong to the spanning tree of $G$, giving at most $3r_k$ possibilities. For each $i$, only the multiset of these $d-2$ choices matters. There are at most $4^m$ choices for the set of cycle edges used, so
\[
    \tiso(G)\leq 4^m\binom{3r_k+d-3}{d-2}^{m}\leq 4^m\left(\frac{C\alpha^k}{d^{5/2}}\right)^{m(d-2)}.
\]
Here we used $r_k=O(\alpha^k k^{-3/2})$ and $\binom{a}{b}\leq(ea/b)^b$. Since $km(d-2)=n-2m$ and $m(d-2)=(1+O(1/d))n/d$, taking logarithms gives the result.
\end{proof}

In the proof of \cref{thm:main}, we made no attempt to optimise the error term. A suitable choice of $k$ and $h$ gives an error of $O(n\log\log d/\log d)$ in the exponent, and allowing the attached rooted trees to have different orders under suitable weighting improves this to $O(n/\log d)$. We conjecture that \cref{prop:necessary-error} gives the correct order of magnitude.

\begin{conjecture}
There is an absolute constant $C>0$ such that every connected $d$-regular graph $G$ on $n$ vertices with $d\geq 2$ satisfies
\[
    \tiso(G)\geq\alpha^n\exp\left(-C\frac{\log d}{d}\,n\right).
\]
\end{conjecture}
It seems plausible that the bound even holds with $C=5/2+o_d(1)$. Indeed, one might expect extremal graphs to resemble our construction, with dense pieces on about $d$ vertices joined by only a few edges.

\subsection*{Acknowledgements.}
The author thanks Julian Sahasrabudhe for his helpful comments on an earlier draft and for his advice on the presentation.

\subsection*{Use of AI}
The proof strategy for the main theorem is due to the author. ChatGPT 6 was used to proofread the manuscript, refine the exposition and review the arguments. It also proposed the construction used to prove \cref{prop:necessary-error}.

\bibliographystyle{plain}
\bibliography{ref}
\end{document}